\documentclass[12pt]{article}
\usepackage{amsmath, latexsym, amssymb}

\numberwithin{equation}{section}

\newcommand{\R}{\mathbb R}
\newcommand{\C}{\mathbb C}

\newcommand{\N}{\mathbb N}

\newcommand{\Bold}[1]{{\boldsymbol{\mathit{#1}}}}

\newcommand{\Z}{\mathbb{Z}}

\newcommand{\QED}{\hspace{.2in}\square\newline}
\newcommand{\qed}{\hspace{.2in}\boxminus\newline}

\newtheorem{theorem}{Theorem}[section]
\newtheorem{corollary}{Corollary}[section]
\newtheorem{proposition}{Proposition}[section]
\newtheorem{definition}{Definition}[section]
\newtheorem{conjecture}{Conjecture}[section]
\newtheorem{lemma}{Lemma}[section]

\begin{document}

\begin{center}
{\Large \textbf{On the infinitude of Prime $k$-tuples}} \vskip 4em

{ J. LaChapelle}\\
\vskip 2em
\end{center}

\begin{abstract}
 Starting with Zhang's theorem  on the infinitude of prime doubles \cite{Z}, we give an inductive argument that there exists an infinite number of prime $k$-tuples for at least one admissible set $\mathcal{H}_k=\{h_1,\ldots,h_k\}$ for each $k$.
\end{abstract}

\section{Introduction}

The eventual historical lesson learned from counting single primes was ``If you want to know the asymptotic behavior of primes, you have to look at $\zeta(s)|_{\Re(s)=1}$". Indeed, it required a thorough understanding of $\zeta(s)$ to finally nail down the prime number theorem (PNT). Even the famous Selberg-Erd\"{o}s `elementary proof' of the PNT conceals $\zeta(s)$ lurking in the background \cite{HA}. Correspondingly, it seems unlikely the Hardy-Littlewood $k$-tuple conjecture \cite{HL} can be settled without possessing a $k$-tuple analog of Riemann zeta.

Nevertheless, number theorists have made impressive gains in the quest to count prime $k$-tuples --- especially recently. As important as these recent advances are, the situation for counting prime $k$-tuples is rather like that for single primes prior to Riemann, Hadamard, and de la Vall\'{e}e Poussin: Without a $k$-tuple zeta function to exploit, the focus has been on showing the infinitude of prime $k$-tuples.

Perhaps the most germane in this respect are Zhang's theorem \cite{Z} and the Maynard-Tao theorem \cite{M},\cite{T} paraphrased by Granville \cite{GR};
\begin{theorem}
There exists an integer $j$ such that the following is true: If
$x+a_1,\ldots,x+a_j$ is an admissible set of forms then there are infinitely many integers
$n$ such that at least two of $n+a_1,\ldots,n+a_j$ are prime numbers.
\end{theorem}
\begin{theorem}
For any given integer $k\geq2$, there exists an integer $j$
such that the following is true: If $x+a_1,\ldots,x+a_j$ is an admissible set of forms then there are infinitely many integers $n$ such that at least $k$ of $n+a_1,\ldots,n+a_j$ are prime numbers.
\end{theorem}

Both of these theorems imply a corollary, again given by Granville \cite{GR};
\begin{corollary}
There is an integer $h$; $0<h\leq B$ such that there are infinitely many
pairs of primes $p,\,p+h$. (for some finite bound $B$)\;.
\end{corollary}
To justify this corollary, choose a sufficiently large finite interval $(0,B]$ and string it together to cover the positive integers in the obvious way. Each interval contains at least two primes out of a finite number of combinations that could occur. Since there are an infinitude of intervals, at least one of those combinations must be represented an infinite number of times. Moreover, such a combination is necessarily admissible because $h$ will be even.

Unfortunately, the same reasoning doesn't work for $k$-tuples with $k>2$. The problem is, one can deduce an infinite number of at least one particular combination but there is no guarantee the combination will constitute an \emph{admissible} $k$-tuple. This is disappointing because belief in the $k$-tuple conjecture is strong, so one strongly expects an infinitude of admissible prime $k$-tuples.

Of course Euclid (and later several others) figured out a way to get the total number of primes without using $\zeta(s)$. And Zhang first did it for (certain)  prime doubles. On the other hand, Euler found a way to utilize $\zeta(s)$ to deduce the infinitude of single primes. Can one generalize Euler to the prime double case and thereby get a handle on a $k$-tuple zeta function? Unfortunately, it is well known the sum over prime-double reciprocals does not diverge. So a straightforward generalization is thwarted from the start.

But maybe a straightforward generalization is not the best approach. Let's briefly re-interpret Euler's method for guidance.
First recall that, if $\mathcal{P}$ is the set of all primes, then
\begin{equation*}
\log(\zeta(s))\asymp\sum_{p\in\mathcal{P}}p^{-s}\;.
\end{equation*}
We aim to show that $\lim_{s\rightarrow1^+}\log(\zeta(s))=\infty$, but suppose we do not know Euler's product representation of $\zeta(s)$. We can still conclude the result by a simple argument.

Assume the contrary. Then $\sum\Lambda(n)/\log(n)n^s$ would converge uniformly at $s=1$, and we would have
\begin{equation*}
\lim_{s\rightarrow1^+}\sum_{n=1}^\infty\frac{\Lambda(n)}{\log(n)n^s}=\lim_{N\rightarrow\infty}\sum_{n\leq N}\frac{\Lambda(n)}{\log(n)n}\;.
\end{equation*}
Hence, partial summation would yield
\begin{equation*}
\lim_{s\rightarrow1^+}\log(\zeta(s))=\log(\log(N_{\mathrm{max}}))-\epsilon+O\left(\frac{1}{\log(N_{\mathrm{max}})}\right)
\end{equation*}
where $\epsilon$ is an end-point contribution and
\begin{equation*}
N_{\mathrm{max}}:=\lim_{N\rightarrow\infty}\sum_{n=1}^N\Lambda(n)=\lim_{N\rightarrow\infty}\sum_{p^k\leq N}\log(p).
\end{equation*}
But from Euclid we know the right-most sum must diverge, and so we have a contradiction.
This $\log(\log(N))$ behavior for $\sum_{p\leq N}1/p$ is well known, and it shows the reciprocal of primes `just barely' diverges with $N$.

The point of reviewing Euler's method is to give a preview of our plan for prime doubles. We don't yet have a representation of the prime-double zeta function $\zeta_{(2)}(s)$. So we will adapt the above argument to the prime-double case and utilize Zhang's result to infer the divergence of a certain pertinent sum for at least one admissible $\mathcal{H}_2=\{0,h\}$. It turns out that the pertinent sum to consider is
\begin{equation*}
\log'(\zeta_{(2)}(s)):=\sum_{n=1}^\infty\frac{\Lambda(n)\Lambda(n+h)}{\log(n)\log(n+h)}\frac{\log(n(n+h)^{1/2})}{(n(n+h))^{s/2}}
\asymp\sum_{p\in\mathcal{P}}\frac{\log(p(p+h))}{(p(p+h))^{s/2}}\;.
\end{equation*}
This sum happens to exhibit the $\log(\log(N))$ behavior.

Of course one could just guess this sum. But it is more satisfying and reassuring to see that it comes from the definition
\begin{equation*}
\log(\zeta_{(2)}(s)):=\sum_{n=1}^\infty\frac{\Lambda(n)\Lambda(n+h)}{\log(n)\log(n+h)}\frac{1}{(n(n+h))^{s/2}}\;,
\end{equation*}
which in turn comes from explicit formulae relating exact and average summatory functions for prime doubles \cite{LA1},\cite{LA2}.

The prime-double case is then extended to $k$-tuples by induction. The reasoning relies crucially on Zhang's theorem: Given an admissible prime $k$-tuple $(p+h_1,\ldots,p+h_k)$, Zhang's theorem implies there exists at least one $h$ such that the amended $(k+1)$-tuple $(p+h_1,\ldots,p+h_k, p+h_k+h)$ is also prime (but not necessarily admissible). Of course this is consistent with the Maynard-Tao theorem.

The final step is to show there must be at least one \emph{admissible} prime $(k+1)$-tuple by this construction. Here we rely on a lemma establishing the fortunate circumstance that $\log'(\zeta_{(2)}(s))\asymp\log'(\zeta_{[2i]}(s))$ where the left side is associated with admissible $h$ and the right side represents an equivalence class of prime doubles $[2i]$ determined by the relation $(2i)^{l'}\sim 2i$ for all $l'\in\mathbb{N}_+$ with $2i$ such that $(2i)^l=h$ for some $i,l\in\mathbb{N}_+$.

The end game and our main result is the following theorem:
\vskip 1em
\noindent\textbf{Theorem} \,\emph{Assume $\mathcal{H}_k=\{0, h_2,\ldots, h_k\}$ is admissible. Then}
\begin{equation}
\lim_{s\rightarrow 1^+}\frac{\sum_{\mathfrak{p}_k\in\mathfrak{P}_k}\frac{\log^{k-1}(p_{(k)})}{p^s_{(k)}}}
{(-1)^{k-1}\log^{(k-1)'}\left(\zeta_{(k)}(s)\right)}=1
\end{equation}
\emph{and}
\begin{equation}\label{divergence}
\lim_{s\rightarrow1^+}(-1)^{k-1}\log^{(k-1)'}\left(\zeta_{(k)}(s)\right)=\infty
\end{equation}
\emph{for at least one admissible $\mathcal{H}_k$.}

\noindent Here $\mathfrak{P}_k$ is the set of admissible prime $k$-tuples, and $p_{(k)}:=(p(p+h_2)\cdots(p+h_k))^{1/k}$ is the geometric mean.

To briefly recap, possessing exact and average summatory functions for admissible prime $k$-tuples allows us to infer certain relevant objects $\log(\zeta_{(k)}(s))$. Although an explicit representation of $\zeta_{(k)}(s)$ remains elusive, $\log(\zeta_{(k)}(s))$ together with Euler's method and Zhang's theorem allow us to deduce an infinitude of at least one admissible prime $k$-tuple for all $k$. History suggests we will need to understand $\zeta_{(k)}(s)$ to go further.

\section{Definitions and lemmas}
This section will introduce some notation/definitions and establish some useful lemmas. But first, some discussion to set the stage.

Studies of patterns in primes often employ sieve methods on sums over the natural numbers $\N_+$ where counting is easy. Our approach instead is to change perspective and consider patterns among prime $k$-tuples on a pair-wise coprime $k$-lattice. The coprime lattice handles {simple} sieving more or less automatically and takes advantage of multiplicative arithmetic functions, but now counting becomes the difficulty. It turns out that the unit counting measure on a ray in the coprime lattice is not the same as the unit measure on $\N_+$: The counting numbers along a coprime ray get renormalized by the density factor $\prod_p(1-{p}^{-1})^{-k}$. Then, further restriction to prime powers brings an additional factor associated with the product of von Mangoldt functions which localize onto prime-power points of the $k$-lattice. Together, these two factors produce the singular series $\mathfrak{S}(\mathcal{H}_k)$ for an admissible prime-power $k$-tuple.\cite{LA2} This supplies a bridge between sums over a ray in the coprime lattice and sums over the natural numbers.

In some sense, one could argue that counting on the coprime lattice is essentially just sieving (at least in spirit). But it turns out that the $k$-lattice perspective brings some new tools to the effort. First of all, as far as coprimes are concerned, ``averages'' are associated with \emph{geometric} means not arithmetic means. More importantly, the coprime lattice hints that prime powers along each dimension are not encoded by  $k$-independent zeta function. Instead, the pattern of explicit sums that one can construct strongly suggests there exists a generalized zeta function $\zeta_{(k)}(s)$ that keeps track of prime power points on the lattice --- in the same way that $\zeta(s)$ does on $\Z_+$.

We are ready for some definitions.

\begin{itemize}
  \item $\mathfrak{P}_k$ is the set of admissible prime $k$-tuples.
  \item $\mathfrak{P}_k\ni\mathfrak{p}_k:=(p,p+h_2,\ldots,p+h_k)$ where $p$ is prime and $\mathcal{H}_k=\{0, h_2,\ldots, h_k\}$ is admissible.
  \item $n_{(k)}:=(n(n+h_2)(n+h_k))^{1/k}$ for integer $n$ is the geometric mean.\footnote{The subscript $(k)$ is supposed to indicate the level $k$ and implicitly the admissible $\mathcal{H}_k$. Sometimes we will make the dependence on $\mathcal{H}_k$ explicit by writing for example $p_{(h_2)}$ for a particular prime double.}
  \item $\mu_{(k)}(n):=(-1)^k\mu(n)\cdots\mu(n+h_k)\,$.
  \item $\lambda_{(k)}(n):=\Lambda(n)\cdots\Lambda(n+h_k)/\log(n)\cdots\log(n+h_k)\,$.
  \item $\log\left(\zeta_{(k)}(s)\right):=\sum_{n=1}^\infty\lambda_{(k)}(n)/n_{(k)}^s\,,\;\;\;\;\Re{(s)}>1\,$.
\end{itemize}

\noindent Now some lemmas.
\begin{lemma}\label{lemma 1}
\begin{equation}\label{bound}
\log\left(\zeta_{(2)}(s)\right)\asymp\sum_{\mathfrak{p}_2\in\mathfrak{P}_2}\frac{1}{p_{(2)}^{s}} =:\log\left(\mathfrak{z}_{(2)}(s)\right)\;,
\end{equation}
and
\begin{equation}
\log'\left(\zeta_{(2)}(s)\right)\asymp\sum_{\mathfrak{p}_2\in\mathfrak{P}_2}\frac{-\log(p_{(2)})}{p_{(2)}^{s}} =\log'\left(\mathfrak{z}_{(2)}(s)\right)\;.
\end{equation}
\end{lemma}
\emph{proof}:
By definition,
\begin{eqnarray}
\log\left(\zeta_{(2)}(s)\right)
&=&\sum_{\mathfrak{p}^\omega_2\in\mathfrak{P}_2}\frac{1}{\omega p^{\omega s/2}}\frac{1}{\omega' (p^\omega+2i)^{\omega's/2}}\;,\;\;\;\;\;\;\Re(s)>1\notag\\
&=&\sum_{\mathfrak{p}_2\in\mathfrak{P}_2}\frac{1}{p^{s/2}}\frac{1}{ (p+2i)^{s/2}}+\sum_{\mathfrak{p}_2\in\mathfrak{P}_2}\sum_{\omega=2}^\infty\sum_{\omega'=2}^\infty\frac{1}{\omega p^{\omega s/2}}\frac{1}{\omega' (p^\omega+2i)^{\omega's/2}}\notag\\\notag\\
&=:&\log\left(\mathfrak{z}_{(2)}(s)\right)+S_{(2)}(s)
\end{eqnarray}
where $\mathfrak{p}^\Bold{\omega}_2:=(p^\omega,(p^\omega+2i)^{\omega'})$ is a prime-power double.
But, for $s=\sigma+it$ with $\sigma, t\in\R$ and $\sigma>1$,
\begin{eqnarray}
\left|S_{(2)}(s)\right|&\leq&\sum_{\mathfrak{p}_2\in\mathfrak{P}_2}\sum_{\omega=2}^\infty\sum_{\omega'=2}^\infty\left|\frac{1}{\omega p^{\omega s/2}}\frac{1}{\omega' (p^\omega+2i)^{\omega's/2}}\right|\notag\\
&<&\sum_{\mathfrak{p}_2\in\mathfrak{P}_2}\sum_{\omega=2}^\infty\sum_{\omega'=2}^\infty\left|\frac{1}{ p^{\omega s/2}}\frac{1}{ (p^\omega+2i)^{\omega's/2}}\right|\notag\\
&<&\sum_{\mathfrak{p}_2\in\mathfrak{P}_2}\sum_{\omega=2}^\infty\sum_{\omega'=2}^\infty\left|\frac{1}{ p^{\omega s/2}}\frac{1}{ (p+2i)^{\omega's/2}}\right|\notag\\
&\leq&\sum_{\mathfrak{p}_2\in\mathfrak{P}_2}\frac{1}{\left(p^\sigma-p^{\sigma/2}\right)}
\frac{1}{\left((p+2i)^\sigma-(p+2i)^{\sigma/2}\right)}\notag\\
&<&\sum_{n=2}^\infty\frac{1}{\left(n^\sigma-n^{\sigma/2}\right)}\frac{1}{\left((n+2i)^\sigma-(n+2i)^{\sigma/2}\right)}\notag\\
&<&\sum_{n=2}^\infty\frac{1}{n^{1/2}(n^{1/2}-1)}\frac{1}{(n+2i)^{1/2}((n+2i)^{1/2}-1)}\notag\\
&<&\sum_{n=2}^\infty\frac{1}{(n-1)^{3/2}}=\zeta(3/2)\;.
\end{eqnarray}

Similarly,
\begin{eqnarray}
S_{(2)}(s)
<\sum_{\mathfrak{p}_2\in\mathfrak{P}_2}\sum_{\omega=2}^\infty\sum_{\omega'=2}^\infty\frac{1}{ p^{\omega s/2}}\frac{1}{ (p+2i)^{\omega's/2}}
&<&\sum_{\mathfrak{p}_2\in\mathfrak{P}_2}\sum_{\omega=2}^\infty\sum_{\omega'=2}^\infty\frac{1}{ (p^{s/2})^{\omega+\omega'}}\notag\\
&=&\sum_{\mathfrak{p}_2\in\mathfrak{P}_2}\frac{1}{p^s\left(p^{s/2}-1\right)^2}\;.
\end{eqnarray}
So, for $M_{(2)}$ a finite positive constant,
\begin{eqnarray}
\left|\frac{dS_{(2)}(s)}{ds}\right|
<\sum_{\mathfrak{p}_2\in\mathfrak{P}_2}\left|\frac{(1-2p^{s/2})\log(p)}{p^s\left(p^{s/2}-1\right)^3}\right|
&<&M_{(2)}\sum_{\mathfrak{p}_2\in\mathfrak{P}_2}\left|\frac{\log(p)}{p^{2s}}\right|\notag\\
&\leq &M_{(2)}\sum_{\mathfrak{p}_2\in\mathfrak{P}_2}\frac{\log(p)}{p^{2\sigma}}\notag\\
&< &M_{(2)}\sum_{n=2}^\infty\frac{\log(n)}{n^{2\sigma}}\notag\\
&< &M_{(2)}\sum_{n=2}^\infty\frac{\log(n)}{n^{2}}\notag\\
&=&-M_{(2)}\,{\zeta}'(2)\;.
\end{eqnarray}

For higher order $k$,
\begin{eqnarray}
S_{(k)}(s)
&<&\sum_{\mathfrak{p}_2\in\mathfrak{P}_2}\sum_{(\omega_1,\ldots,\omega_k)=2}^\infty
\frac{1}{ p^{\omega_1 s/2}}\frac{1}{ (p+h_2)^{\omega_2s/2}}\cdots\frac{1}{ (p+h_k)^{\omega_ks/2}}\notag\\
&<&\sum_{\mathfrak{p}_2\in\mathfrak{P}_2}\sum_{(\omega_1,\ldots,\omega_k)=2}^\infty
\frac{1}{ (p^{s/2})^{\omega_1+\cdots+\omega_k}}\notag\\
&=&\sum_{\mathfrak{p}_2\in\mathfrak{P}_2}\frac{1}{p^{ks/2}\left(p^{s/2}-1\right)^k}\;,
\end{eqnarray}
and
\begin{eqnarray}
\left|\frac{d^{(k-1)}S_{(k)}(s)}{ds^{(k-1)}}\right|
<M_{(k)}\sum_{\mathfrak{p}_2\in\mathfrak{P}_2}\left|\frac{\log(p)^{k-1}}{p^{ks}}\right|
&\leq &M_{(k)}\sum_{\mathfrak{p}_2\in\mathfrak{P}_2}\frac{\log(p)^{k-1}}{p^{k\sigma}}\notag\\
&< &M_{(k)}\sum_{n=2}^\infty\frac{\log(n)^{k-1}}{n^{k\sigma}}\notag\\
&< &M_{(k)}\sum_{n=2}^\infty\frac{\log(n)^{k-1}}{n^{k}}\notag\\
&=&(-1)^{k-1}M_{(k)}\,\zeta^{(k-1)'}(k)
\end{eqnarray}
with $M_{(k)}$ bounded.
$\qed$

The following two lemmas utilize a particularly useful interpretation of sums of the form $\sum \mu_{(2)}(n)\Lambda_{(2)}(n)f(n)$. Since the pre-factor localizes onto prime doubles \cite{LA1}, it is advantageous to view the sum as a double sum over the coprime $2$-lattice.

\begin{lemma}\label{lemma 2}
Denote the coprime $2$-lattice by $\{(n_1,n_2)\in\mathbb{N}_+^2 \;|\; \mathrm{gcd}(n_1,n_2)=1\}$, and let $2i<2j=(2i)^l$ with fixed $l\in\mathbb{N}_+$. Then
\begin{equation}
\sum_{n=1}^{\infty}\frac{\mu_{(2i)}(n)\lambda_{(2i)}(n)}{n^{s/2}(n+2i)^{s/2}}\log(n(n+2i))
-\sum_{n=1}^{\infty}\frac{\mu_{(2j)}(n)\lambda_{(2j)}(n)}{n^{s/2}(n+2j)^{s/2}}\log(n(n+2j))\,<\,(2j-2i)\;.
\end{equation}

\end{lemma}
\emph{proof}: For $N>1+(2j-2i)$,
\begin{eqnarray}\label{compare}
\sum_{n=1}^{N}\frac{\mu_{(2i)}(n)\lambda_{(2i)}(n)}{n^{s/2}(n+2i)^{s/2}}\log(n(n+2i))\notag\\
&&\hspace{-2.75in}\propto\sum_{n_1= 1}^{N}\sum_{n_2= 1+2i}^{N+2i}
\frac{\mu(n_1)\lambda_{(1)}(n_1)\mu(n_2)\lambda_{(1)}(n_2)}{n_1^{s/2}n_2^{s/2}}\log(n_1n_2)
\,\delta(n_2\,,\, n_1+2i)\notag\\
&&\hspace{-2.75in}=\sum_{n_1=1}^{N}\,\frac{\mu(n_1)\lambda(n_1)}{n_1^{s/2}}
\left(\sum_{n_2= 1+2i}^{N+2i}
\frac{\mu(n_2)\lambda(n_2)}{n_2^{s/2}}\,\log(n_1n_2)\,\delta(n_2\,,\, n_1+2i)\right)\notag\\
&&\hspace{-2.75in}=\sum_{n_1=1}^{N}\,\frac{\mu(n_1)\lambda(n_1)}{n_1^{s/2}}
\left\{\left(\sum_{n_2= 1+2j}^{N+2j}
\frac{\mu(n_2)\lambda(n_2)}{n_2^{s/2}}\,\log(n_1n_2)\,\delta(n_2\,,\, n_1+2i)\right)
+r(n_1;2i,2j,N)\right\}\notag\\
&&\hspace{-2.75in}=\sum_{n_1=1}^{N}\,\frac{\mu(n_1)\lambda(n_1)}{n_1^{s/2}}
\left\{\left(\sum_{n_2= 1+2j}^{N+2j}
\frac{\mu(n_2)\lambda(n_2)}{n_2^{s/2}}\,\log(n_1n_2)\,\delta(n_2\,,\, n_1+2j)\right)
+r(n_1;2i,2j,N)\right\}\notag\\
&&\hspace{-2.75in}\propto
\sum_{n=1}^{N}\frac{\mu_{(2j)}(n)\lambda_{(2j)}(n)}{n^{s/2}(n+2j)^{s/2}}\log(n(n+2j))+R(2i,2j,N)
\end{eqnarray}
where
\begin{equation}
r(n_1;2i,2j,N):=
\left(\sum_{n_2=1+2i}^{2j}-\sum_{n_2=N+1+2i}^{N+2j}\right)
\frac{\mu(n_2)\lambda(n_2)}{n_2^{s/2}}\,\log(n_1n_2)\,\delta(n_2\,,\, n_1+2i)
\end{equation}
and
\begin{equation}
R(2i,2j,N):=
\sum_{n_1=1}^{N}\,\frac{\mu(n_1)\lambda(n_1)}{n_1^{s/2}}r(n_1;2i,2j,N)\;.
\end{equation}

To bound $R(2i,2j,N)$ interchange summation order to get
\begin{eqnarray}
|R(2i,2j,N)|&\leq&\left|\sum_{p_2=1+2i}^{2j}\frac{\log(p(p-2i))}{2(p(p-2i))^{s/2}}\right|
+\left|\sum_{p_2=N+1+2i}^{N+2j}\frac{\log(p(p-2i))}{2(p(p-2i))^{s/2}}\right|\;,\;\;\;\;\Re(s)>1\notag\\
&<&\sum_{n_2=1+2i}^{2j}\left(\frac{1}{2}+\frac{1}{2}\right)\notag\\
&=& 2j-2i\;.
\end{eqnarray}

Returning now to (\ref{compare}), note the delta function in the second line restricts the double sum to the appropriate ray $\mathbf{R}_{(2i)}$ in the pair-wise coprime $2$-lattice. According to \cite{LA2}, the proportionality constant is $\prod_p(1-p^{-1})^{2}$. The final line transforms back to $\N_+$ from the coprime lattice picking up the inverse proportionality constant $\prod_p(1-p^{-1})^{-2}$. The result follows as $N\rightarrow\infty$ since both series converge for $\Re{(s)}>1$ and the proportionality constants cancel.
$\qed$

By the same token,
\begin{lemma}\label{lemma 3}
Let $2j<2m=(2j)^l$ with fixed $l\in\mathbb{N}_+$, then
\begin{equation}
\sum_{n=1}^{\infty}\frac{\mu_{(2j)}(n)\lambda_{(2j)}(n)}{n^{s/2}(n+2j)^{s/2}}\log(n(n+2j))
-\sum_{n=1}^{\infty}\frac{\mu_{(2m)}(n)\lambda_{(2m)}(n)}{n^{s/2}(n+2m)^{s/2}}\log(n(n+2m))
\,<\,2m-2j\;.
\end{equation}
\end{lemma}
\emph{proof}:
The proof follows the same argument as the preceding lemma.
$\qed$

Therefore choosing either  $2i<2j=(2i)^l$ or $2i=(2j)^l> 2j$,
\begin{eqnarray}
\left|\log'\left(\mathfrak{z}_{(2i)}(s)\right)-\log'\left(\mathfrak{z}_{(2j)}(s)\right)\right|
&<&|2j-2i|
\end{eqnarray}
follows from Lemma \ref{lemma 2} and Lemma \ref{lemma 3}.
Hence,
\begin{equation}
\log'\left(\mathfrak{z}_{[2i]}(s)\right)\asymp\log'\left(\mathfrak{z}_{(2j)}(s)\right),\;\;\;\;\Re{(s)}>1
\end{equation}
where the equivalence class includes all integer powers of $2i$. Together with Lemma \ref{lemma 1} conclude
\begin{corollary}\label{equivalence}
\begin{equation}
\log'\left(\zeta_{[2i])}(s)\right)\asymp\log'\left(\zeta_{(2j)}(s)\right),\;\;\;\;\Re{(s)}>1\;.
\end{equation}
\end{corollary}
This result makes sense because the associated sums are along congruent rays in the pair-wise coprime lattice, and (heuristically at least) there are factors of $C_{[2i])}$ and $C_{(2j)}$ coming from the singular series which are asymptotically equivalent since $2j=(2i)^l$. 

\section{Euler's lead}
With these preliminaries, we can apply Euler's method to prime doubles. 

Recall the definition
\begin{definition}
\begin{equation}
\log\left(\zeta_{(2)}(s)\right):=\sum_{n=1}^\infty\frac{\lambda_{(2)}(n)}{n_{(2)}^s}\,,\;\;\;\;\Re{(s)}>1\;.
\end{equation}
\end{definition}
We have
\begin{proposition}\label{proposition}
Let $2j=h$ where $h$ comes from Zhang's theorem. Then
\begin{equation}
\lim_{s\rightarrow1^+}\log'\left(\zeta_{(2j)}(s)\right)=\infty\;.
\end{equation}
\end{proposition}
\emph{Proof}:
Assume the contrary. Then $\lim_{s\rightarrow1^+}\log'\left(\zeta_{2j}(s)\right)$ is bounded and so converges. Consequently, for every $\varepsilon>0$ there exists an $M$ such that $m>M$ implies
\begin{equation}
\left|\sum_{n=m+1}^{m+l}\frac{\lambda_{(2j)}(n)}{(n+2j)}\log(n_{(2j)})\right|<\varepsilon
\end{equation}
for each $l\in\mathbb{N}_+$. It follows that $\log'\left(\zeta_{2j}(s)\right)$ converges uniformly for $\Re{(s)}\geq 1$.

Now, the PNT implies $\lambda_{(1)}(n+2j)=\Lambda(n+2j)/\log(n+2j)=O(1/\log(n+2j))$. Moreover, $\log(n_{(2)})/\log(n)= O(1)$. So by uniform convergence, the PNT, and partial summation we get
\begin{eqnarray}
\lim_{s\rightarrow1^+}\log'\left(\zeta_{(2j)}(s)\right)
&>&-\lim_{N\rightarrow\infty}\sum_{n=1}^N\frac{\lambda_{(2j)}(n)}{(n+2j)}\log(n_{(2j)})\notag\\
&\sim&-\lim_{N\rightarrow\infty}
\sum_{n=1}^N\frac{w_{(2)}(n)\Lambda(n+2j)}{(n+2j)\log(n+2j)}\notag\\
&=&\log(\log(N_{\mathrm{max}}^{(j)}))-\epsilon_{2j}+O\left(\frac{1}{\log(N_{\mathrm{max}}^{(j)})}\right)
\end{eqnarray}
where the weight  $w_{(2)}(n)=1$ if $\Lambda(n)\neq0\wedge\Lambda(n+2j)\neq0$ and $w_{(2)}(n)=0$ otherwise, $N_{\mathrm{max}}^{(j)}:=\lim_{N\rightarrow\infty}\sum_{n=1}^{N}w_{(2)}(n)\Lambda(n+2j)$, and the constant $\epsilon_{2j}$ is an inconsequential end-point contribution. But the PNT and Zhang's theorem imply $N_{\mathrm{max}}^{(j)}=\infty$, and we arrive at a contradiction.
$\QED$

Therefore, for at least one $2j$, Proposition \ref{proposition} and Corollary \ref{equivalence} imply
\begin{corollary}
\begin{equation}
\lim_{s\rightarrow1^+}\log'\left(\mathfrak{z}_{[2i]}(s)\right)
=\sum_{\mathfrak{p}_2\in\mathfrak{P}_2}\frac{\log(p_{[2i]})}{p_{[2i]}}
=\infty
\end{equation}
where $[2i]$ is the equivalence class determined by the relation $(2i)^{l'}\sim 2i$ for all $l'\in\mathbb{N}_+$ with $2i$ such that $(2i)^l=2j$ for some $i,l\in\mathbb{N}_+$.
\end{corollary}
The sum must include an infinite number of terms, and so there are infinitely many $\mathfrak{p}_{2}\in\mathfrak{P}_2$ for each admissible $\mathcal{H}_{2}=\{0,[2i]\}$.

\section{Proof of Theorem}
We restate the theorem for easy reference.
\begin{theorem}\label{theorem}
Assume $\mathcal{H}_k$ is admissible. Then
\begin{equation}
\lim_{s\rightarrow 1^+}\frac{\sum_{\mathfrak{p}_k\in\mathfrak{P}_k}\frac{\log^{k-1}(p_{(k)})}{p^s_{(k)}}}
{(-1)^{k-1}\log^{(k-1)'}\left(\zeta_{(k)}(s)\right)}=1
\end{equation}
and
\begin{equation}\label{divergence}
\lim_{s\rightarrow 1^+}(-1)^{k-1}\log^{(k-1)'}\left(\zeta_{(k)}(s)\right)=\infty
\end{equation}
for at least one admissible $\mathcal{H}_k$.
\end{theorem}

\emph{Proof}: It is straightforward to show that (\ref{bound}) generalizes to $\log\left(\zeta_{(k)}(s)\right)\asymp\sum_{\mathfrak{p}_k\in\mathfrak{P}_k}p_{(k)}^{-s}$
with $\left|S(s)\right|<\zeta(\frac{k+1}{k})$ and recall $|S^{(k-1)'}(s)|$ is bounded. So
\begin{equation}
\sum_{\mathfrak{p}_k\in\mathfrak{P}_k}\frac{\log^{k-1}(p_{(k)})}{p_{(k)}^{s}}
\asymp(-1)^{k-1}\log^{(k-1)'}\left(\zeta_{(k)}(s)\right)\;.
\end{equation}

Now let $w_{(k+1)}(n)=1$ encode the condition that $\Lambda_{(k)}(n)\neq0\wedge\Lambda(n+h_{k+1})\neq0$ and $w_{(k+1)}(n)=0$ otherwise. Note that $w_{(k+1)}(n)= w_{(k)}(n)w_{(1)}(n+h_{k+1})$. 
Also note that
\begin{eqnarray}
\log^{(k-1)'}\left(\zeta_{(k)}(s)\right)
&=&(-1)^{k-1}\sum_{n=1}^\infty\frac{\lambda_{(k)}(n)}{n^s_{(k)}}\log^{k-1}(n_{(k)})\notag\\
&>&(-1)^{k-1}\sum_{n=1}^\infty\frac{\lambda_{(k)}(n)}{(n+h_{k})^s}\log^{k-1}(n_{(k)})\;.
\end{eqnarray}
Evidently, assuming the first $k-1$ elements of $(n,n+h_2,\ldots,n+h_k)$ are all prime and that $\lambda_{(k)}(n)=O(\lambda_{(1)}(n+h_k)\log^{1-k}(n_{(k)}))$  for $\mathcal{H}_k$ admissible is tantamount to assuming (\ref{divergence}) which in turn implies infinitely many $\mathfrak{p}_k\in\mathfrak{P}_k$.

So, adopting the assumption that $\lambda_{(k+!)}(n)=O(\lambda_{(1)}(n+h_{k+1})\log^{-k}(n_{(k+1)}))$ given that $w_{(k)}(n)=1$ for some admissible $\mathcal{H}_k$, get
\begin{eqnarray}\label{condition}
(-1)^k\log^{(k)'}\left(\zeta_{(k+1)}(s)\right)
&=&\sum_{n=1}^\infty\frac{\lambda_{(k+1)}(n)}{n^s_{(k+1)}}
\log^k(n_{(k+1)})\notag\\
&>&\sum_{n=1}^\infty\frac{\lambda_{(k+1)}(n)}{(n+h_{k+1})^s}
\log^k(n_{(k+1)})\notag\\
&\sim&\sum_{n=1}^\infty
\frac{w_{(k+1)}(n)\Lambda(n+h_{k+1})}{\log(n+h_{k+1})(n+h_{k+1})^s}\notag\\
\end{eqnarray}
where the last line uses the assumption as well as $\lambda_{(1)}(n+h_{k+1})=O(1/\log (n+h_{k+1}))$ implied by the PNT.

As in the prime-double case, assume $(-1)^k\log^{(k)'}\left(\zeta_{(k+1)}(s)\right)$ is bounded at $s=1$. Partial summation yields
\begin{equation}
\lim_{s\rightarrow 1^+}(-1)^k\log^{(k)'}\left(\zeta_{(k+1)}(s)\right)
\sim\log(\log(\mathrm{max}\,{ x}))-\epsilon_{h_{k+1}}\;.
\end{equation}
We again have a contradiction since $\mathrm{max}\,{ x}=\lim_{N\rightarrow\infty}\sum_{n=1}^{N}w_{(k+1)}(n)
\Lambda(n+h_{k+1})=\infty$ follows from: (i) the assumption $w_{(k)}(n)=1$, (ii) $w_{(k+1)}(n)= w_{(k)}(n)w_{(1)}(n+h_{k+1})$, and (iii) Zhang's work for at least one $2j_{k+1}$ with $0<2j_{k+1}\leq h_{k+1}-h_k$ provided $h_{k+1}-h_k$ is chosen large enough. Hence, $\lim_{s\rightarrow 1^+}(-1)^k\log^{(k)'}\left(\zeta_{(k+1)}(s)\right)=\infty$ for the associated equivalence class $[2i_{k+1}]$ by previous arguments.

If it happens that $2i_{k+1}\equiv0$ (mod $k+1$), then $h_k+[2i_{k+1}]$ and $h_k$ belong to the same residue class  (mod $k+1$) so $\mathcal{H}_{k+1}=\{0,\ldots,h_k,h_k+[2i_{k+1}]\}$ is automatically admissible. Conversely, if $2i_{k+1}\equiv\!\!\!\!\not \;\;0$ (mod $k+1$), then $\mathcal{H}_{k+1}$ can be rendered admissible for a suitable choice of representative in $[2i_{k+1}]$. For example, if $2i_{k+1}\equiv\!\!\!\!\not \;\;0$ (mod $k+1$), then $h_{k+1}=h_k+(2i_{k+1})$ and $h_{k+1}=h_k+(2i_{k+1})^2$ belong to different residue classes  (mod $k+1$). Consequently at least one of them will yield an admissible $\mathcal{H}_{k+1}$, because at least one case will not occupy the complete set of residue classes modulo primes.

Finally, Zhang's result guarantees the induction assumption is true at $k=2$ for at least one $\mathcal{H}_2=\{0,2j_2\}$. It follows that $\lim_{s\rightarrow 1^+}(-1)^k\log^{(k-1)'}\left(\zeta_{(k)}(s)\right)=\infty$ for at least one admissible $\mathcal{H}_{k}=\{0,[2i_2],[2i_2]+[2i_3],\ldots,\sum_{\hat{k}\leq k}[2i_{\hat{k}}]\}$  for all $k$ by induction.
$\QED$

\begin{corollary}
Given an admissible $\mathcal{H}_k$ satisfying Theorem \ref{theorem}, if some $2j_{\hat{k}}=2^a$ for integer $a$, then there are infinitely many twin primes.
\end{corollary}

\emph{Proof}: The $(0,\ldots,\hat{k},\ldots)$ components of $\mathfrak{p}_k$ determine an admissible $\mathcal{H}_2=\{0,2^a\}$ and its associated prime doubles $\mathfrak{p}_2$ of which there are necessarily infinitely many. Then previous arguments imply the claim. Of course similar statements for other types of prime doubles can be made for any component of $\mathfrak{p}_k$ with $2j_{\hat{k}}=(2i_{\hat{k}})^a$.
$\QED$

One would like to extend this theorem to \emph{all} admissible $\mathcal{H}_k$. But such an extension would require one to show that $\mathrm{max}\,{ x}\rightarrow\infty$ for \emph{all} prime doubles. In this case, the assumption used for the induction argument could be made for all admissible $\mathcal{H}_k$ and verified at $k=2$. Using ideas from \cite{LA2} leads to

\begin{conjecture}
If $\zeta_{(2)}(s)$ is meromorphic on $\C$ with a simple pole at $s=1$ and no zeros on $\Re(s)=1$, then Theorem \ref{theorem} holds for all admissible $\mathcal{H}_k$.
\end{conjecture}


\begin{thebibliography}{99}

\bibitem{Z}
Yitang Zhang, Bounded gaps between primes, \emph{Ann. Math.}, \textbf{179}(3), (2014), 1121--1174.

\bibitem{HA}
D. Harvey, Selberg's Symmetyry Formula, \emph{Expo. Math.}, \textbf{22}, (2004), 185-195.

\bibitem{HL}
G.H Hardy and J.E. Littlewood, Some problems of `Partitio numerorum'; III: On the expression of a number as a sum of primes, \emph{Acta Math.}, \textbf{44}(1), (1923), 1--70.


\bibitem{M}
J. Maynard, Small Gaps Between Primes, arXiv:math.NT/1311.4600v2 (2013).

\bibitem{T}
T. Tao, Polymath8b: Bounded intervals with many primes, after Maynard,
http://terrytao.wordpress.com/2013/11/19/polymath8b-bounded-intervals-with-many-primes-after-maynard/

\bibitem{GR} A. Granville, Primes in Intervals of Bounded Length,
www.dms.umontreal.ca/~andrew/CEBBrochureFinal.pdf

\bibitem{LA1}
J. LaChapelle, Exact Summatory Functions for Prime $k$-tuples, arXiv:math.NT/1406.5533 (2014).

\bibitem{LA2}
J. LaChapelle, A Gamma Distribution Hypothesis for Prime $k$-tuples, arXiv:math.NT/1406.6289 (2014).




\end{thebibliography}
\end{document}